\documentclass[oneside,english]{amsart}
\usepackage[T1]{fontenc}
\usepackage[utf8]{inputenc}
\usepackage{amstext}
\usepackage{amsthm}
\usepackage{amssymb}

\makeatletter
\numberwithin{equation}{section}
\numberwithin{figure}{section}
\theoremstyle{plain}
\newtheorem{thm}{\protect\theoremname}
\theoremstyle{plain}
\newtheorem{lem}[thm]{\protect\lemmaname}
\theoremstyle{plain}
\newtheorem{prop}[thm]{\protect\propositionname}

\makeatother

\usepackage{babel}
\providecommand{\lemmaname}{Lemma}
\providecommand{\propositionname}{Proposition}
\providecommand{\theoremname}{Theorem}

\begin{document}
\title{A note on Barker Sequences of even length}
\author{Jürgen Willms}
\email{willms.juergen@fh-swf.de}
\address{Institut für Computer Science, Vision and Computational Intelligence,
Fachhochschule Südwestfalen, Lindenstraße 53, D-59872 Meschede, Germany}
\begin{abstract}
A Barker sequence is a binary sequence for which all nontrivial aperiodic
autocorrelations are either 0, 1 or -1. The only known Barker sequences
have length 2, 3, 4, 5, 7, 11 or 13. It is an old conjecture that
no longer Barker sequences exist and in fact, there is an overwhelming
evidence for this conjecture. For binary sequences of odd length,
this conjecture is known to be true, whereas for even length it is
still open, whether a Barker sequence of even length greater 4 exists. 

Similar to the well-known fact that a Barker sequence of odd length
is necessarily skew-symmetric, we show that in the case of even length
there is also a form of symmetry albeit weaker. In order to exploit
this symmetry, we derive different formulas for the calculation of
the aperiodic correlation. We prove by using only elementary methods
that there is no Barker sequence of even length $n>4$ with $C_{1}=C_{3}=\cdots=C_{\frac{n}{2}-1},$
where $C_{k}$ denotes the $k$th aperiodic autocorrelation of the
sequence. 
\end{abstract}

\keywords{binary sequence, Barker sequence, autocorrelation}
\subjclass[2000]{11B83, 05B10, 94A55}
\date{1.4.2021}

\maketitle

\section{Introduction\label{sec:Introduction}}

For an integer $n>1$ an element of $\{-1,1\}^{n}$ is called a\emph{
binary sequence }$a$ of \emph{length} $n.$ We write $a_{j}$ for
the $j$-th entry in $a$. For $k=0,1,\cdots,n-1$ the $k$th \emph{aperiodic
autocorrelation} of the binary sequence $a$ is defined by
\[
C_{k}(a)=\sum_{j=1}^{n-k}a_{j}a_{j+k}.
\]

The autocorrelation C$_{0}(a)$ is called the\emph{ peak autocorrelation}
of $a$, whereas all other values $C_{k}(a)$ are called the \emph{nontrivial}
aperiodic autocorrelations. Note that $C_{0}(n)=n$. In the following
$a$ will always denote a binary sequence of length $n>1$ and we
will abbreviate $C_{k}(a)$ as $C_{k}$. A simple parity argument
shows that 
\begin{equation}
C_{k}\equiv n-k\pmod2\label{eq:simpleParity}
\end{equation}

for each $k=1,2,\cdots,n-1.$ Furthermore, we have
\begin{equation}
(\sum_{j=1}^{n}a_{j})^{2}=n+2\sum_{k=1}^{n-1}C_{k}.\label{eq:sum(a)^2}
\end{equation}

Moreover, use the simple parity argument $xy\equiv x-y+1\pmod4$ for
$x,y\in\{-1,1\}$ to obtain for $k=1,2,\cdots,n-1$
\begin{equation}
C_{k}+C_{n-k}=\sum_{j=1}^{n}a_{j}a_{j+k}\equiv n\pmod4,\label{eqn:periodicMod4}
\end{equation}
where the index $j+k$ of $a$ is reduced modulo $n$ if $j+k>n$.

Of particular interest are binary sequences for which the absolute
values of all aperiodic autocorrelations are small; for a recent survey
see for example \cite{schmidt2016sequences}. An ideal sequence from
this point of view is a Barker sequence. A binary sequence $a$ is
called \emph{Barker sequence} if all nontrivial aperiodic autocorrelations
of $a$ are either 0, 1 or -1. Hence by (\ref{eq:simpleParity}) all
nontrivial aperiodic autocorrelations of a Barker sequence are in
magnitude as small as possible. 

The only known Barker sequences have length 2, 3, 4, 5, 7, 11 or 13.
It is a longstanding conjecture that no longer Barker sequences exist.
If $n$ is odd, then a well-known result of Turyn and Storer \cite{turyn1961binary}
published in 1961 says, that Barker sequences of odd length n > 13
do not exist. More then 50 years later, it was noted \cite{willms2014counterexamples}
that their elementary but arduous inductive proof is incomplete. A
correction of the proof was suggested in \cite{schmidt2016barker}.
A new proof was given by Borwein and Erdélyi \cite{borwein2013note}
in 2013 and an even simpler and very different proof was presented
by Schmidt and Willms \cite{schmidt2016barker} in 2016. 

The Barker conjecture for binary sequences of even length is still
open, but there is an overwhelming evidence for it. For example, results
from Borwein and Mossinghoff \cite{borwein2014wieferich} and Leung
and B. Schmidt \cite{leung2016anti} show that there is no Barker
sequence of even length $n$ for $4<n\leq4\cdot10^{33}.$

Let us briefly try to explain why in the Barker conjecture the case
of even length appears to be considerably harder than the case of
odd length. Suppose that $a$ is a Barker sequence of length $n$
and let $1\leq k<n$. If the length $n$ is odd then by (\ref{eq:simpleParity})
we know that $C_{k}=0$ if $k$ is odd, and (\ref{eqn:periodicMod4})
gives us that $C_{k}=(-1)^{\frac{n-1}{2}}$ if $k$ is even; so in
this case all aperiodic autocorrelations are fixed. In contrast, if
the length $n$ of a Barker sequence is even, then $C_{k}=0$ if $k$
is even. But now (\ref{eqn:periodicMod4}) gives us only that
\begin{equation}
C_{k}=-C_{n-k}.\label{eq:Cor_antisymmetric}
\end{equation}
So in this case only half of the aperiodic autocorrelations that are
not zero are fixed; no further constraints on the other half (as for
example $C_{1},C_{3},\cdots,C_{\frac{n}{2}-1}$) are known (apart
of course from the fact that each value is either 1 or -1). 

Furthermore, by (\ref{eq:sum(a)^2}) if $a$ is a Barker sequence
of even length $n\geqslant4$ then $|\sum_{j=1}^{n}a_{j}|=\sqrt{n}$
and hence $n=4r^{2}$ for some $r\in\mathbb{N}$; in particular, $n$
is a multiple of 4. There are many more stronger and deeper partial
results on the non-existence of Barker sequences of even length (see
for example \cite{jedwab2008can} or \cite{borwein2014wieferich}),
which however we will not use in the following. 

Let $a$ be a Barker sequences of even length $n\geq4$ and suppose
that the first $\frac{n}{4}$ odd aperiodic autocorrelations $C_{1},C_{3},\cdots,C_{\frac{n}{2}-1}$
are all equal. The next result shows that this is only possible if
the length $n$ of $a$ equals 4. 
\begin{thm}
\label{thm:Ck_positive}Suppose $a$ is a Barker sequences of even
length $n\geq4$ with $C_{1}=C_{3}=\cdots=C_{\frac{n}{2}-1}$. Then
$n=4$. 
\end{thm}

We need a couple of lemmata in order to prove Theorem \ref{thm:Ck_positive}
and proceed as follows. In the Section 2, we consider what we call
weak symmetric sequences, and show that a putative Barker sequence
of even length $n\geq4$ must be weak symmetric. In Section 3, we
derive different formulas for the calculation of the aperiodic correlation
that are in particular useful in the weak symmetric case. Finally,
in Section 4 we apply our results to sequences with a stronger form
of symmetry, which enables us to prove Theorem \ref{thm:Ck_positive}. 

Let us outline the main idea behind the proof of Theorem \ref{thm:Ck_positive}.
Suppose that $a$ is a putative Barker sequence of even length $n\geq4$
with $C_{1}=C_{3}=\cdots=C_{\frac{n}{2}-1}$ and put $u=\frac{n}{4}.$
Similar as in \cite{turyn1961binary} we show that 
\begin{equation}
a_{j}a_{j+1}=a_{2j}a_{2j+1}\text{ for each \ensuremath{j} satisfying \ensuremath{1\leq j<u.}}\label{eq:zweier-1}
\end{equation}
Now construct a new sequence $p$ from $a$ by reversing the first
half of $a$ and do the same with the second half of $a$. It follows
from Lemma \ref{lem:C(k)_b_reverse} that $p$ is also a Barker sequence
of length $n$. Applying again (\ref{eq:zweier-1}) gives us that
$|C_{u}(b)|=u\neq0$ where $b$ denotes the binary sequence of length
$2u$ with $b_{j}=a_{j}$ for each $1\leq j\leq2u.$ Since by Lemma
\ref{lem:C(k)_strong} we have $C_{k}(b)=0$ for even $k$, it follows
that $u$ must be odd.\footnote{The conclusion that $u$ is odd, here derived by merely elementary
arguments, follows also from a result of Turyn \cite{turyn1965character}
using fairly deep methods of character theory.} Now assume that $n>4.$ Then it is not difficult to show that both
$C_{2}(b)$ and $C_{2u-2}(b)$ must be zero. But by (\ref{eqn:periodicMod4})
we have $C_{2}(b)+C_{2u-2}(b)\equiv2u\pmod4$ which implies that $u$
is even, a contradiction. Hence $n=4.$ 

\section{Weak Symmetry}

In the following, we will exploit the symmetry of a putative Barker
sequence of even length. A very similar approach was used in the case
of odd length Barker sequences; all proofs in this case rely heavily
on the fact that a Barker sequence of odd length $n$ is skew-symmetric,
that is $a_{j}a_{n+1-j}=(-1)^{\frac{n+1}{2}+j}$ for each $1\leq j\leq n.$
In the case of a Barker sequence of even length we will prove a weaker
form of symmetry that leads us to the following definition. We call
a binary sequence $a$ \emph{weak symmetric} if the length $n$ of
$a$ is a multiple of 4 and $a_{j}a_{j+1}=-a_{n+1-j}a_{n-j}$ holds
for each odd integer $j$ with $1\leq j<\frac{n}{2}$ (and thus for
each odd integer $j$ with $1\leq j<n).$ In Proposition \ref{prop:barker1}
we will show that a Barker sequence of even length $n\geq4$ is indeed
weak symmetric.

For $k=1,2,\cdots,n$ put

\[
T_{k}=\sum_{j=1}^{k}a_{j}a_{n+1-j}.
\]

The next result shows that weak symmetry, $T_{k}$ and $C_{n-k}$
are strongly related.
\begin{lem}
\label{lem:T(k)}Let $a$ be a binary sequence of length $n$ and
let \textup{$1\leq k<n$. Then:}

\begin{enumerate}

\item $C_{n-k}\equiv T_{k}\pmod4.$

\item  Suppose that a is weak symmetric. Then $T_{k}=0$ if $k$
is even, and $T_{k}=a_{k}a_{n+1-k}$ if $k$ is odd. 

\end{enumerate}
\end{lem}

\begin{proof}
(i) Use the simple counting argument that is stated in Lemma \ref{lem:Counting}
(i) to obtain

\[
(-1)^{\frac{k-T_{k}}{2}}=\prod_{j=1}^{k}a_{j}a_{n+1-j}=\prod_{j=1}^{k}a_{j}a_{j+n-k}=(-1)^{\frac{k-C_{n-k}}{2}}
\]
which implies that $C_{n-k}\equiv T_{k}\pmod4$.

(ii) Suppose $a$ is weak symmetric. Then $a_{1}a_{2}=-a_{n-1}a_{n}$
and thus $T_{2}=0.$ Since for $k>2$ we have $T_{k}=T_{k-2}+a_{k-1}a_{n-k}+a_{k}a_{n+1-k}$,
it follows that $T_{k}=0$ if $k$ is even. Set $T_{0}=0;$ if $k$
is odd, then $k=2j+1$ for some $j\geq0$ and thus $T_{k}=T_{2j}+a_{k}a_{n+1-k}=a_{k}a_{n+1-k}.$
\end{proof}
The next lemma uses a simple counting argument. Recall that $b_{j}$
denotes the $j$-th entry in $b$.
\begin{lem}
\label{lem:Counting} Let $b$ be a binary sequence of length $m$.

\begin{enumerate}

\item Let $S=\sum_{j=1}^{m}b_{j}$; then
\[
\prod_{j=1}^{m}b_{j}=(-1)^{\frac{m-S}{2}}
\]

\item  Let $R=\sum_{j=1}^{m}(-1)^{j}b_{j}$; then
\[
\prod_{j=1}^{m}b_{j}=\begin{cases}
(-1)^{\frac{R}{2}} & \text{if \ensuremath{m} is even }\\
(-1)^{\frac{R+1}{2}} & \text{\text{if \ensuremath{m} is odd}}.
\end{cases}
\]

\end{enumerate}
\end{lem}

\begin{proof}
(i) Denote by $\alpha$ resp. $\beta$ the number of positive resp.
negative entries in $b$. Then $\alpha+\beta=m$, $S=\alpha-\beta$
and hence $m-S=2\beta$. Since $\prod_{j=1}^{m}b_{j}=(-1)$$^{\beta}$,
the assertion (i) follows directly.

(ii) Let $\gamma$ denote the cardinality of the set $\left\{ 1\leq j\leq m:j\text{ odd}\right\} $.
Then $\gamma=\frac{m}{2}$ if $m$ is even, and $\gamma=\frac{m+1}{2}$
if $m$ is odd. Hence we have $(-1)^{\gamma}\prod_{j=1}^{m}b_{j}=\prod_{j=1}^{m}(-1)^{j}b_{j}$.
Use (i) to obtain $\prod_{j=1}^{m}(-1)^{j}b_{j}=(-1)^{\frac{m-R}{2}}.$
By considering the cases $m$ even and $m$ odd separately, assertion
(ii) follows easily.
\end{proof}
We have already seen in Section \ref{sec:Introduction} that the length
$n$ of a Barker sequence of even length is is a multiple of 4 and
that its nontrivial aperiodic autocorrelations $C_{k}$ equal zero
if $k$ is even; moreover, it is also weak symmetric as the next result
shows.
\begin{prop}
\label{prop:barker1}Let $a$ be a Barker sequence of even length
$n\geq4.$\textup{ Then }$a$ is weak symmetric, and $C_{k}=-a_{k}a_{n+1-k}$
if $k$ is odd with \textup{$1\leq k<n.$} 
\end{prop}

\begin{proof}
Let $1\leq j\leq\frac{n}{2}$. Note that $T_{1}=a_{1}a_{n}=C_{n-1}$
and $|T_{j+1}-T_{j}|=1.$ Suppose that $a$ is a Barker sequence of
even length $n\geq4$. Then $|C_{n-j}-C_{n-j-1}|.$ Since by Lemma
\ref{lem:T(k)} (i) we have $C_{n-j}\equiv T_{j}\pmod4$, it follows
that $C_{n-j}=T_{j}.$ In particular, we have $T_{j}=0$ if $j$ is
even. If $j$ is odd then $0=T_{j+1}-T_{j-1}=a_{j}a_{n+1-j}+a_{j+1}a_{n-j},$
which shows that $a$ is weak symmetric. If $k$ is odd with $1\leq k<n$
then by Lemma \ref{lem:T(k)} (ii) and (\ref{eq:Cor_antisymmetric})
we have $C_{k}=-C_{n-k}=-T_{k}=-a_{k}a_{n+1-k}.$
\end{proof}

\section{Correlation Formulas }

In order to simplify the notation and to exploit the symmetry we will
use the following abbreviation. For a binary sequence $a$ of length
$n$ put $\delta_{k}=a_{k}a_{n+1-k}$ for $k=1,2,\cdots,n.$

Note that for $1\leq k\leq n$ we have
\begin{equation}
\delta_{k}=\delta_{n+1-k}.\label{eq:delta}
\end{equation}

Moreover, if the binary sequence $a$ of length $n$ is weak symmetric,
then it follows directly from the definition that for each odd $k$
(with $1\leq k<n$) 
\begin{equation}
\delta_{k}=-\delta_{k+1}.\label{eq:delta_weak_sym}
\end{equation}

One goal of this section is to express the aperiodic autocorrelations
$C_{k}$ of a binary sequence $a$ of length $n=2m$ by using only
the first $m$ entries $a_{1},a_{2},\cdots,a_{m}$ of $a$ together
with $\delta_{1},\delta_{2},\cdots,\delta_{m}$. The next lemma is
a first step towards this goal. 
\begin{lem}
\label{lem:C(k)_simple} Let $a$ be a binary sequence $a$ of length
$n$ and let \textup{$1\leq k<n$}. 

\begin{enumerate}

\item If $k$ is even then

\[
C_{n-k}=\sum_{j=1}^{\frac{k}{2}}a_{j}a_{k+1-j}(\delta_{j}+\delta_{k+1-j}).
\]

\item  If $k$ is odd then 

\[
C_{n-k}=\delta_{\frac{k+1}{2}}+\sum_{j=1}^{\frac{k-1}{2}}a_{j}a_{k+1-j}(\delta_{j}+\delta_{k+1-j}).
\]

\end{enumerate}
\end{lem}

\begin{proof}
By definition we have $C_{n-k}=\sum_{j=1}^{k}a_{j}a_{j+n-k}=\sum_{j=1}^{k}a_{j}a_{k+1-j}\delta_{k+1-j}.$
In order to prove (i) and (ii) combine in the last sum the first with
the last summand, the second with the second last, and so on.
\end{proof}
Let $a$ be a binary sequence $a$ of even length $n=2m$ and suppose
that $1\leq k\leq m$. Note that the formula for $C_{n-k}$ in Lemma
\ref{lem:C(k)_simple} does not refer to an element $a_{j}$ or $\delta_{j}$
with $j>m.$ In order to to achieve the same goal when computing $C_{k}$
(instead of $C_{n-k}$) we have to split the sum on the right hand
side of (i) resp. (ii) of Lemma \ref{lem:C(k)_simple} into two parts
as shown in the next lemma.
\begin{lem}
\label{lem:C(k)_composite} Let $a$ be a binary sequence $a$ of
even length $n=2m$ and let \textup{$1\leq k<m$}. 

\begin{enumerate}

\item If $k$ is even then

\[
C_{k}=\sum_{j=1}^{m-k}a_{j}a_{j+k}(1+\delta_{j}\delta_{j+k})+\sum_{j=1}^{\frac{k}{2}}a_{j+m-k}a_{m+1-j}(\delta_{j+m-k}+\delta_{m+1-j}).
\]

\item  If $k$ is odd then 

\[
C_{k}=\delta_{\frac{n-k+1}{2}}+\sum_{j=1}^{m-k}a_{j}a_{j+k}(1+\delta_{j}\delta_{j+k})+\sum_{j=1}^{\frac{k-1}{2}}a_{j+m-k}a_{m+1-j}(\delta_{j+m-k}+\delta_{m+1-j}).
\]

\end{enumerate}
\end{lem}

\begin{proof}
We show that both (i) and (ii) follow directly from Lemma \ref{lem:C(k)_simple}.
Suppose first that $k$ is even. Then Lemma \ref{lem:C(k)_simple}
(i) gives us that $C_{k}=\sum_{j=1}^{\frac{n-k}{2}}a_{j}a_{n-k+1-j}(\delta_{j}+\delta_{n-k+1-j}).$
If $1\leq j\leq m-k$ then by (\ref{eq:delta}) we have $a_{n-k+1-j}=a_{k+j}\delta_{k+j},$
$\delta_{k+j}=\delta_{n-k+1-j},$ and hence $a_{j}a_{n-k+1-j}(\delta_{j}+\delta_{n-k+1-j})=a_{j}a_{j+k}\delta_{k+j}(\delta_{j}+\delta_{n-k+1-j})=a_{j}a_{j+k}(\delta_{j}\delta_{j+k}+1).$
Since 
\[
\sum_{j=m-k+1}^{\frac{n-k}{2}}a_{j}a_{n-k+1-j}(\delta_{j}+\delta_{n-k+1-j})=\sum_{j=1}^{\frac{k}{2}}a_{j+m-k}a_{m+1-j}(\delta_{j+m-k}+\delta_{m+1-j})
\]
assertion (i) follows. Suppose next that $k$ is odd. Then by Lemma
\ref{lem:C(k)_simple} (ii) we have $C_{k}=\delta_{\frac{n-k+1}{2}}+\sum_{j=1}^{\frac{n-k-1}{2}}a_{j}a_{n-k+1-j}(\delta_{j}+\delta_{n-k+1-j}),$
and the same argument as above can be used to prove (ii).
\end{proof}

\section{\label{sec:Strong-Symmetric }Strong Symmetric Binary Sequences}

In the this section we will use Lemma \ref{lem:C(k)_simple} and Lemma
\ref{lem:C(k)_composite} in order to derive a formula for the aperiodic
correlation if the binary sequence exhibits a stronger form of symmetry,
which is for example the case if $a$ is a Barker sequence of even
length $n=2m$ $\geq4$ with $C_{1}=C_{3}=\cdots=C_{m-1}=1.$ 

Suppose that $a$ is a Barker sequence of even length $n=2m\geq4$
with $C_{1}=C_{3}=\cdots=C_{m-1}=1$ and let $1\leq k\leq m$. By
Proposition \ref{prop:barker1} we have $a_{k}a_{n+1-k}=-1$ if $k$
is odd, and $a_{k}a_{n+1-k}=1$ if $k$ is even. This motivates the
following definition.

A binary sequence $a$ of length $n$ is called \emph{strong symmetric}
if the length $n$ of $a$ is a multiple of 4 and and $a_{j}a_{n+1-j}=(-1)^{j}$
for each $1\leq j\leq\frac{n}{2}$.

Note that if $a$ is a strong symmetric binary sequence of length
$n=2m$ then $a$ is weak symmetric and $\delta_{k}=(-1)^{k}$ for
each $1\leq k\leq m$, and $\delta_{k}=(-1)^{k+1}$ for each $m<k\leq n$
(recall that $\delta_{k}$ is defined as $\delta_{k}=a_{k}a_{n+1-k}$
for $1\leq k\leq n$). Note further, that a strong symmetric binary
sequence $a$ of length $n=2m$ is already determined by the first
$m$ entries $a$$_{1}$$,a_{2},\cdots,a_{m}$.

Applying Lemma \ref{lem:C(k)_simple} and Lemma \ref{lem:C(k)_composite}
to strong symmetric binary sequences gives us the following result. 
\begin{lem}
\label{lem:C(k)_strong} Let $a$ be a strong symmetric binary sequence
of length $n=2m$. Denote by $b$ be the binary sequence of length
$m$ with $b_{j}=a_{j}$ for each \textup{$1\leq j\leq m.$ }If $k$
is even with \textup{$1\leq k<n,$} then 

\[
C_{k}=\begin{cases}
2C_{k}(b) & \text{if \ensuremath{k<m} }\\
0 & \text{\text{if \ensuremath{k\geq m}}}.
\end{cases}
\]

Furthermore, if $k$ is odd with $1\leq k<m$, then following two
statements hold:

\begin{enumerate}

\item 
\[
C_{n-k}=(-1)^{\frac{k+1}{2}}+2\sum_{j=1}^{\frac{k-1}{2}}(-1)^{j}b_{j}b_{k+1-j}.
\]

\item 
\[
C_{k}=(-1)^{\frac{k-1}{2}}-2\sum_{j=1}^{\frac{k-1}{2}}(-1)^{j}b_{j+m-k}b_{m+1-j}.
\]

\end{enumerate}
\end{lem}

Note that $C_{k}(b)$ denotes the $k$th aperiodic autocorrelation
of the binary sequence $b$, whereas $C_{k}$ as always refers to
the sequence $a$, i.e. $C_{k}=C_{k}(a)$. Note further that the sequence
$b$ corresponds to the first half of the sequence $a$, so equations
(i) and (ii) hold also, if we replace on the right hand side the entries
of $b$ with the corresponding entries of $a$.
\begin{proof}
Let $1\leq k<n$. Suppose first that $k$ is even. If $k\geq m$ then
$C_{k}=0$ follows directly from Lemma \ref{lem:C(k)_simple} (i),
since in this case $\delta_{j}+\delta_{n-k+1-j}=0$ for each $1\leq j\leq\frac{n-k}{2}$.
For the case $k<m$ we use Lemma \ref{lem:C(k)_composite} (i). Note
that in this case $1+\delta_{j}\delta_{j+k}=2$ for each $1\leq j\leq m-k$
and $\delta_{j+m-k}+\delta_{m+1-j}=0$ for each $1\leq j\leq\frac{k}{2}$;
hence $C_{k}=2\sum_{j=1}^{m-k}a_{j}a_{j+k}=2C_{k}(b).$ 

Let us suppose next that $k$ is odd and $1\leq k<m$. Then by Lemma
\ref{lem:C(k)_simple} (ii) we have $C_{n-k}=\delta_{\frac{k+1}{2}}+\sum_{j=1}^{\frac{k-1}{2}}a_{j}a_{k+1-j}(\delta_{j}+\delta_{k+1-j}).$
Note that in this case $\delta_{j}+\delta_{k+1-j}=2(-1)^{j}$ for
all $1\leq j\leq\frac{k-1}{2}$ and that $\delta_{\frac{k+1}{2}}=(-1)^{\frac{k+1}{2}}$,
which proves (i). In order to prove (ii) use Lemma \ref{lem:C(k)_composite};
for this note that $\delta_{\frac{n-k+1}{2}}=(-1)^{\frac{k-1}{2}}$,
$1+\delta_{j}\delta_{j+k}=0$ for each $1\leq j\leq m-k$, and since
$m$ is even we have $\delta_{j+m-k}+\delta_{m+1-j}=-2(-1)^{j}$ for
each $1\leq j\leq\frac{k-1}{2}.$ 
\end{proof}
Suppose that $a$ and $p$ are both strong symmetric binary sequences
of length $n=2m$. If the first $m$ entries of $p$ are the same
as the first $m$ entries of $a$ but in reversed order, then the
aperiodic autocorrelations of $a$ and $p$ are strongly related;
this shows the next lemma.
\begin{lem}
\label{lem:C(k)_b_reverse} Let $a$ and $p$ be strong symmetric
binary sequences of length $n=2m$ and suppose that $p_{j}=a_{m+1-j}$
for each $1\leq j\leq m$. Then for each $1\leq k<n$ \textup{
\[
C_{k}(p)=\begin{cases}
\qquad C_{k}(a) & \text{if \ensuremath{k} is even }\\
-C_{n-k}(a) & \text{\text{if \ensuremath{k} is odd}}.
\end{cases}
\]
}
\end{lem}

\begin{proof}
Let $1\leq k<n,$ and as in Lemma \ref{lem:C(k)_strong} let $b$
denote the binary sequence of length $m$ with $b_{j}=a_{j}$ for
each $1\leq j\leq m.$ Consider first the case when $k\text{ is even.}$
Then by Lemma \ref{lem:C(k)_strong} we have $C_{k}(a)=0=C_{k}(p)$
if $k\geq m.$ Since reversing a binary sequence does not chance its
aperiodic autocorrelations, it follows that $C_{k}(a)=2C_{k}(b)=C_{k}(p)$
if $k<m$.

Next consider the case when $k\text{ is odd.}$ If $k<m$ then by
Lemma \ref{lem:C(k)_strong} (ii) and (i) we have $C_{k}(p)=(-1)^{\frac{k-1}{2}}-2\sum_{j=1}^{\frac{k-1}{2}}(-1)^{j}a_{k+1-j}a_{j}$
and $C_{n-k}(a)=(-1)^{\frac{k+1}{2}}+2\sum_{j=1}^{\frac{k-1}{2}}(-1)^{j}a_{j}a_{k+1-j}.$
This shows that $C_{k}(p)=-C_{n-k}(a)$ if $k<m$. The same argument
gives us also that $C_{k}(a)=-C_{n-k}(p)$ if $k<m.$ Hence $C_{k}(p)=-C_{n-k}(a)$
for all odd $k$ satisfying $1\leq k\leq n$.
\end{proof}
Recall that if $a$ is a Barker sequence of even length $n\geq4$
then $n$ is a multiple of 4. We now consider the special case of
a Barker sequence of even length $n=4u$ where the first $u$ odd
aperiodic autocorrelations all have the same value.
\begin{lem}
\label{lem:2er_regel} Suppose that $a$ is a Barker sequence of even
length $n=2m\geq4$ with $C_{1}=C_{3}=\cdots=C_{m-1},$ and let \textup{$u=\frac{n}{4}$}.
Then we have $a_{j}a_{j+1}=a_{2j}a_{2j+1}=a_{u+j}a_{u+j+1}$ for each
$j$ satisfying $1\leq j<u$. 
\end{lem}

\begin{proof}
Recall that by Proposition \ref{prop:barker1} the Barker sequence
$a$ is weak symmetric. Consider first the case that $C_{1}=C_{3}=\cdots=C_{m-1}=1$,
and let $1\leq k<m.$ Then by Proposition \ref{prop:barker1} we have
$\delta_{k}=a_{k}a_{n+1-k}=-1$ if $k$ is odd; since $a$ is weak
symmetric this implies that $a$ is in fact strong symmetric. If $k$
is odd then by (\ref{eq:Cor_antisymmetric}) and Lemma \ref{lem:C(k)_strong}
(i) we have $-1=C_{n-k}=(-1)^{\frac{k+1}{2}}+2R_{k}$ with $R_{k}=\sum_{j=1}^{\frac{k-1}{2}}(-1)^{j}a_{j}a_{k+1-j}.$
Depending on whether $\frac{k+1}{2}$ is odd or even, $R_{k}$ is
either $0$ or $-1.$ In both cases Lemma \ref{lem:Counting} (ii)
gives us that $\prod_{j=1}^{\frac{k-1}{2}}a_{j}a_{k+1-j}=1$ and hence
$\prod_{j=1}^{k}a_{j}=a_{\frac{k+1}{2}}$ for each odd $k$ with 1$\leq k<m.$
Similar as in \cite{turyn1961binary} multiplying two successive equations
of this latter form shows that 
\begin{equation}
a_{j}a_{j+1}=a_{2j}a_{2j+1}\text{ for each \ensuremath{j} satisfying \ensuremath{1\leq j<u.}}\label{eq:zweier}
\end{equation}

Next consider the case $C_{1}=C_{3}=\cdots=C_{m-1}=-1.$ Since the
binary sequence $s$ of length $n$ defined by $s_{j}=(-1)^{j}a$$_{j}$
is a Barker sequence with $C_{1}(s)=C_{3}(s)=\cdots=C_{m-1}(s)=1$
Equation (\ref{eq:zweier}) holds also in this case. 

It remains to show that $a_{2j}a_{2j+1}=a_{u+j}a_{u+j+1}$ for $1\leq j<u$.
As in Lemma \ref{lem:C(k)_b_reverse} let $p$ be the (well defined)
strong symmetric binary sequence of length $n$ with $p_{j}=a_{m+1-j}$
for each $1\leq j\leq m.$ Lemma \ref{lem:C(k)_b_reverse} shows that
$p$ is also a Barker sequence of length $n$ with $C_{1}(p)=C_{3}(p)=\cdots=C_{m-1}(p)=-C_{1}$.
Thus by (\ref{eq:zweier}) we have $a_{m+1-k}a_{m-k}=a_{m+1-2k}a_{m-2k}$
for $1\leq k<u$. Putting $j=u-k$ this can be rewritten as $a_{u+j+1}a_{u+j}=a_{2j+1}a_{2j}$
for $1\leq j<u$. 
\end{proof}
We are now in the position to prove Theorem \ref{thm:Ck_positive}.
\begin{proof}[Proof of Theorem \ref{thm:Ck_positive}]
 Suppose that $a$ is a Barker sequence of even length $n=2m\geq4$
with $C_{1}=C_{3}=\cdots=C_{m-1}.$ Without loss of generality we
can assume that $C_{1}=1.$ Otherwise, instead of $a$ consider the
binary sequence $s$ of length $n$ defined by $s_{j}=(-1)^{j}a$$_{j}$;
then, as mentioned before, $s$ is also a Barker sequence of length
$n$ with $C_{1}(s)=C_{3}(s)=\cdots=C_{m-1}(s)=1.$ 

Put $u=\frac{n}{4}$ and as in Lemma \ref{lem:C(k)_strong} let $b$
denote the binary sequence of length $m$ with $b_{j}=a_{j}$ for
each $1\leq j\leq m.$ Then by Lemma \ref{lem:2er_regel} we have
$a_{j}a_{j+u}=a_{j+1}a_{j+1+u}$ for all $j$ satisfying $1\leq j<u$
and hence $|C_{u}(b)|=u\neq0.$ Thus by Lemma \ref{lem:C(k)_strong}
the integer $u$ must be odd. 

Now assume that $m>2.$ Then Lemma \ref{lem:C(k)_strong} gives us
that $C_{2}(b)=0$ and $C_{m-2}(b)=0$. But by (\ref{eqn:periodicMod4})
we have $0=C_{2}(b)+C_{m-2}(b)\equiv m\pmod4$ which implies that
$u$ is even, a contradiction. Hence $m=2.$ 
\end{proof}

\section{Conclusion}

We have shown that Barker sequences of even length $n\geq4$ exhibit
certain symmetry properties. In order to exploit this symmetry, different
formulas for the calculation of the aperiodic correlation have been
derived. By using only elementary methods we have proven that the
special case where $C_{1}=C_{3}=\cdots=C_{\frac{n}{2}-1}$ is for
Barker sequences of even length $n\geq4$ only possible if $n=4$.
It is an interesting open question whether a similar elementary approach
can be used to further reduce the number of cases for which the Barker
conjecture for even length can be proved.


\bigskip{}

\end{document}